\documentclass[12pt,reqno]{amsart}
\usepackage{amscd,amsmath,amsopn,amssymb,amsthm,multicol}
\usepackage{hyperref}
\DeclareMathOperator{\Ad}{Ad}
\DeclareMathOperator{\ad}{ad}

\newcommand{\thickline}{\noalign{\hrule height 1pt}}

\newtheorem{theorem}{Theorem}
\newtheorem{lemma}{Lemma}

\newtheorem{proposition}{Proposition}
\newtheorem{defidition}{Definition}
\newtheorem{example}{Example}
\newtheorem{corollary}{Corollary}

\begin{document}

\title[ Homogeneous Geodesics in Generalized Wallach Spaces]
{ Homogeneous Geodesics in Generalized Wallach Spaces}
\author{Andreas Arvanitoyeorgos and Yu Wang}

\address{University of Patras, Department of Mathematics, GR-26500 Rion, Greece}
\email{arvanito@math.upatras.gr}
\address{Sichuan university of Science and Engineering, Zigong, 643000, China}
\email{wangyu\_813@163.com}

\begin{abstract}
We classify  generalized Wallach spaces which are g.o. spaces. We also investigate homogeneous geodesics in generalized Wallach spaces for any given invariant Riemannian metric and we give some examples.

 \medskip


\medskip

 \end{abstract}
\maketitle

\section*{Introduction}
 Let $(M, g)$  be a homogeneous Riemannian manifold, i.e. a connected Riemannian manifold on which the largest connected group $G$ of isometries acts transitively. Then $M$ can be expressed as a homogeneous space $(G/K, g)$ where $K$ is the isotropy group at a fixed pointed $o$ of $M$, and $g$ is a $G$-invariant metric. In this case the Lie algebra $\mathfrak{g}$ of $G$ has an $\Ad(K)$-invariant  decomposition $\mathfrak{g}=\mathfrak{k}\oplus\mathfrak{m}$, where $\mathfrak{m} \subset \mathfrak{g}$ is a linear subspace of $\mathfrak{g}$ and $\mathfrak{k}$  is the Lie algebra of $K$. In general such decomposition is not  unique. The $\ad(K)$-invariant subspace $\mathfrak{m}$ can be naturally identified with the tangent space $T_{o}M$ via the projection $\pi:G\rightarrow G/K$.

 A geodesic $\gamma(t)$ through the origin $o$ of $M=G/K$ is called {\it homogeneous} if it is an orbit of a one-parameter subgroup of $G$, that is
\begin{equation}\label{1}
 \gamma(t)=\mathrm{exp}(tX)(o), \quad t \in \mathbb{R},
 \end{equation}
where $X$ is a non zero vector of $\mathfrak{g}$.

A homogeneous Riemannian manifold is called a {\it g.o. space}, if all geodesics are homogeneous with respect to the largest connected group of isometries. All naturally reductive spaces  are g.o. spaces (\cite{Ko-No}), but the converse is not true in general. In \cite{Kap} A. Kaplan proved the existence of g.o. spaces that are in no way naturally reductive.  These are generalized Heisenberg groups with two-dimensional center.
In \cite{Ko-Pr-Va} O. Kowalski, F. Pr\" ufer and L. Vanhecke made an explicit classification of all naturally reductive spaces up to dimension five. In \cite{Ko-Va} O. Kowalski and L. Vanhecke gave a classification of all g.o. spaces, which are in no way naturally reductive, up to dimension six.
In \cite{Gor} C. Gordon described g.o. spaces  which are nilmanifolds and in \cite{Tam} H. Tamaru  classified homogeneous g.o. spaces which are fibered over irreducible symmetric spaces.
In  \cite{Du2} and \cite{Du-Ko1} O. Kowalski and Z. Du\v sek investigated  homogeneous geodesics in   Heisenberh groups  and some $H$-type groups. Examples of g.o. spaces in dimension seven were obtained by Du\v sek, O. Kowalski and S. Nik\v cevi\' c in (\cite{Du-Ko-Ni}). Also, in \cite{Al-Ar} the first  author and D.V. Alekseevsky  classified generalized flag manifolds which are   g.o. spaces.

 Concerning the existence of homogeneous geodesics in  homogeneous Riemannian manifold, we recall the following.
 V.V. Kajzer  proved that a Lie group endowed with a left-invariant metric admits at least one homogeneous geodesic (\cite{Kaj}). O. Kowalski and J. Szenthe extended this result to all homogeneous Riemannian manifolds (\cite{Ko-Sz}). An extension of the result of \cite{Ko-Sz} to reductive homogeneous  pseudo-Riemannian manifolds has been also obtained (\cite{Du-Ko1}, \cite{Phi}).
 Also, O. Kowalski, S. Nik\v cevi\'c and Z. Vl\' a\v sek studied homogeneous geodesics in homogneous Riemannian manifolds (\cite{Ko-Ni-Vl}), and  G. Calvaruso  and R. Marinosci studied homogeneous geodesics in three-dimension Lie groups (\cite{Mar},  \cite{Ca-Ma}).
 Homogeneous geodesics were studied by J. Szenthe (\cite{Sze1}, \cite{Sze2}, \cite{Sze3}, \cite{Sze4}).
 In addition, D. Latifi studied homogeneous geodesics in homogeneous Finsler spaces (\cite{Lat}), and the first author  investigated homogeneous  geodesics in  the flag manifold $SO(2l+1)/U(l-m)\times SO(2m+1)$ (\cite{Arv}).

Homogeneous geodesics in the affine setting were studied in \cite{Du-Ko-Vl}.
Finally, D.V. Alekseevsky and Yu. G. Nikonorov in \cite{Al-Ni} studied the structure of compact g.o. spaces and gave some sufficient conditions
for existence and non existence of an invariant metric  with homogeneous geodesics on
a homogeneous space of a compact Lie group G. They also gave a classification of compact simply
connected g.o. spaces  of positive Euler characteristic.

 Because of these results, it is natural to study g.o. spaces as well as to describe  homogeneous geodesics for other large classes of homogeneous spaces. In this paper, we study this problem for generalized Wallach spaces.  These spaces were well known before as three-locally-symmetric spaces (\cite{Lo-Ni-Fi}), however they were recently  classified by Z. Chen, Y. Kang and K. Liang (\cite{Ch-Ka-Li}) and Yu.G. Nikonorov (\cite{Nik3}).
 We search for homogeneous geodesics in these spaces.

 One of the main results in the present paper is Theorem \ref{T3}, which classifies generalized Wallach spaces which are g.o. spaces.
 For those which are not, we show how to obtain homogeneous geodesics. We make explicit computations for the three dimensional Lie group $SU(2)$ (thus recovering a result of R.A. Marinosci), and for the Stiefel manifold $SO(4)/SO(2)$.

 The paper is organised as follows:  In Section 1 we  recall the basic definitions and properties of homogeneous geodesics in a Riemannian manifold. In Section 2 we  recall the definition of generalized Wallach spaces as well as their classification from \cite{Nik3}. In Section 3 we   classify g.o. spaces among generalized Wallach spaces.  For those which are not g.o. spaces, in Section 4 we discuss how to find all homogeneous geodesics  for a given $G$-invariant Riemannian metric.

 \medskip

\noindent
{\bf Acknowledgements.}  The first author was supported by Grant \# E.037 from the research committee of the
University of Patras.

\noindent
The second author is supported by NSFC 11501390,  by  Sichuan Province University
Key Laboratory of Bridge Non-destruction Detecting and Engineering Computing 2014QZJ03,  by funding of the Education
Department of Sichuan Province 15ZB0218, 15ZA0230, 15SB0122, and by funding of Sichuan University of Science and Engineering grant
2014PY06, 2015RC10.

\noindent
The second author also acknowledges the hospitality of the University of Patras, during her visit 2015-2016.

\section{Homogeneous geodesics in homogeneous Riemannian manifolds}\label{Section1}

 Let $(M=G/K, g)$ be a homogeneous Riemannian manifold. Let $\mathfrak{g}$ and $\mathfrak{k}$ be the Lie algebras of $G$ and $K$ respectively and let
\begin{equation}\label{2}
  \mathfrak{g}=\mathfrak{k}\oplus \mathfrak{m}
\end{equation}
be a reductive decomposition. The canonical projection $\pi:G \rightarrow G/K $ induces an isomorphism between the subspace $\mathfrak{m}$ and the tangent space $T_{o}M$ at the identity $o=eK$. The $G$-invariant metric $g$   induces a scalar product $\langle \cdot,\cdot \rangle$ on $\mathfrak{m}$ which is $\Ad({K})$-invariant. Let $B(\cdot, \cdot)=-$Killing form on $\mathfrak{g}$. Then any $\Ad({K})$-invariant scalar product $\langle \cdot,\cdot \rangle$ on $\mathfrak{m}$ can be expressed as $\langle x, y \rangle=B(\Lambda x, y)\; (x, y \in \mathfrak{m})$, where $\Lambda$ is an $\Ad({K})$-equivariant positive definite symmetric operator on $\mathfrak{m}$. Conversely, any such operator $\Lambda$ determines an  $\Ad({K})$-invariant scalar product $\langle x, y \rangle=B(\Lambda x, y)$ on $\mathfrak{m}$, which in turn determines   a $G$-invariant Riemannian metric $g$ on $\mathfrak{m}$. We say that  $\Lambda$ is the {\it operator associated} to the metric $g$, or simply the {\it associated operator}.
Also, a Riemannian metric generated by inner
  product $B(\cdot,\cdot)$ is called {\it standard metric}.

\begin{defidition} \label {D1}  A nonzero vector $X \in \mathfrak{g}$ is called a geodesic vector if the curve (\ref{1}) is a geodesic.

\end{defidition}

\begin{lemma}[\cite{Ko-Va}]  \label{L1} A nonzero vector  $X \in \mathfrak{g}$ is a geodesic vector if and only if

\begin{equation}\label{3}
  \langle [X,Y]_{\mathfrak{m}},X_{\mathfrak{m}} \rangle =0,
\end{equation}

for all $Y \in \mathfrak{m}$.
Here the subscript $\mathfrak{m}$ denotes the projection into $\mathfrak{m}$.
\end{lemma}

A useful description of homogeneous geodesics (\ref{1}) is provided by the following :
\begin{proposition} [\cite{Al-Ar}]  \label{P1}

 Let $(M=G/K, g)$ be a homogeneous Riemannian manifold and $\Lambda$ be the associated operator. Let $a\in \mathfrak{k}$ and $x \in \mathfrak{m}$. Then the following are equivalent:

 (1)\  The orbit $\gamma(t)=\mathrm{exp}t(a+x)\cdot o$ of the one-parameter subgroup $\mathrm{exp}t(a+x)$ through the point $o=eK$ is a geodesic of $M$.

 (2) \ $[a+x, \Lambda x] \in \mathfrak{k}$.

 (3) \ $\langle [a, x], y \rangle = \langle x, [x, y]_{\mathfrak{m}} \rangle \ \mbox {for all} \  y \in \mathfrak{m}$.

 (4) \ $\langle [a+x, y]_{\mathfrak{m}}, x \rangle=0 \  \mbox {for all} \ y \in \mathfrak{m}$.
 \end{proposition}

 An important corollary of Proposition \ref{P1} is the following:

 \begin{corollary} [\cite{Al-Ar}]  \label{C1}

 Let $(M=G/K, g)$ be a homogeneous Riemannian manifold. Then $(M=G/K, g)$ is a g.o. space if and only if for every $x\in \mathfrak{m}$  there exists an $a(x) \in \mathfrak{k}$  such that
 $$
 [a(x)+x, \Lambda x] \in \mathfrak{k}.
 $$
 \end{corollary}

  If we want to decide  whether a  homogeneous Riemannian manifold $(M=G/K, g)$ is a  g.o. space or not, we need to   find a decomposition of the form (\ref{2}) and look for geodesic vectors of  the form
  \begin{equation}\label{4}
  X=\sum^{s}_{i=1}a_{i}e_{i}+\sum^{l}_{j=1}x_{j}A_{j}.
  \end{equation}
  Here $\{e_{i}: i=1,2, \dots, s\}$ is a  basis of $\mathfrak{m}$ and $\{A_{j}: j=1,2,\dots,l\}$ is a  basis of $\mathfrak{k}$.
  By substituting $X$ into equation (\ref{3}) we obtain a system of linear algebraic equations for the parameters $x_j$. If for any  $X \in \mathfrak{m}\setminus \{0\}$  this system for the  variables $x_j$ has real solutions, then it follows that  the  homogeneous Riemannian manifold $(M=G/K, g)$ is  a g.o. space.

  If, on the other hand, we need to find all homogeneous geodesics in the homogeneous Riemannian manifold $(M=G/K, g)$, then we have to calculate all geodesic vectors in the Lie algebra $\mathfrak{g}$. Condition (\ref{3}) reduces to  a system of $s$   quadratic equations for the variables $x_j$ and $a_i$.  Then the geodesic vectors correspond to those solutions of  this system for the variables $x_1,\dots, x_l,  a_1, \dots, a_s$,  which are not all equal to zero.

 \section{Generalized Wallach spaces }

Let $G/K$ be a compact homogeneous space with connected compact semisimple Lie group $G$ and a compact subgroup $K$. Denoted by $\mathfrak{g}$ and  $\mathfrak{k}$ Lie algebras of $G$ and $K$ respectively. We assume that $G/K$ is almost effective, i.e., there are no non-trial ideals of the Lie algebra
  $\mathfrak{g}$ in $\mathfrak{k}\subset \mathfrak{g}$. Let  $\langle \cdot,\cdot \rangle=B(\Lambda\cdot,\cdot)$ be an $\Ad({K})$-invariant scalar product on $\mathfrak{m}$, where  $\Lambda$ is the associated operator.

  We assume that the homogeneous space $G/K$ has the following property. The module $\mathfrak{m}$  decomposes into a direct sum of three $\Ad(K)$-invariant irreducible modules pairwise orthogonal with respect to $B$, i.e.
  \begin{equation}\label{5}
  \mathfrak{m}=\mathfrak{m}_1\oplus\mathfrak{m}_2\oplus\mathfrak{m}_3,
  \end{equation}
  such that

 \begin{equation}\label{6}
  [\mathfrak{m}_i, \mathfrak{m}_i]\subset\mathfrak{k} \quad  i=1,2,3.
  \end{equation}
  A homogeneous space with this property are called {\it generalized Wallach space}.\\

  Some examples of these spaces are  the manifolds of completely flags in the complex, quaternionic, and Cayley projective planes, that is
  $SU(3)/{T_{max}}$, $Sp(3)/Sp(1)\times Sp(1) \times Sp(1)$,  $F_4/Spin(8)$ (known as Wallach spaces), and the generalized flag manifolds  $SU(n_1+n_2+n_3)/S(U(n_1)\times U(n_2)\times U(n_3))$, $SO(2n)/U(1)\times U(n-1)$ and $\quad E_6/U(1)\times U(1)\times SO(8)$.
  
   Every generalized Wallach space admits a three parameter family of invariant Riemannian metrics determined by $\Ad({K})$-invariant inner products
  \begin{equation}\label{7}
  \langle \cdot,\cdot\rangle=\lambda_1B(\cdot,\cdot)\mid_{\mathfrak{m}_1}+\lambda_2B(\cdot,\cdot)\mid_{\mathfrak{m}_2}+\lambda_3B(\cdot,\cdot )\mid_{\mathfrak{m}_3},
  \end{equation}
  where $\lambda_1, \lambda_2, \lambda_3$ are positive real numbers.

  Let $d_i$ be the dimension of $\mathfrak{m}_i$. Let $\{e^j_{i}\}$ be an orthonormal basis of  $\mathfrak{m}_j$ with respect to $B$, where $j=1,2,3\ $ and $ 1 \leq i \leq d_j$. Consider the expression $[ijk]$ defined by the equality
  \begin{equation}\label{8}
  [ijk]=\sum_{\alpha,\beta,\gamma}B( [e^{i}_{\alpha},e^{j}_{\beta}],e^{k}_{\gamma})^2,
  \end{equation}
  where $\alpha,\beta,\gamma$ range from 1 to $ d_i, d_j $ and $d_k$ respectively (cf. \cite{Wa-Zi}). The symbols $[ijk]$ are symmetric in all three indices due to  the bi-invariance  of the metric $\langle \cdot,\cdot \rangle$. Moreover, for the  spaces under consideration we have $[ijk]=0$,  if two induces coincide.

  We recall the classification of generalized Wallach spaces that was recently obtained by Yu.G. Nikoronov (\cite{Nik3}) and  Z. Chen, Y. Kang, K. Liang (\cite{Ch-Ka-Li}):

\begin{theorem} [\cite{Nik3}, \cite{Ch-Ka-Li}]\label{T1} Let $G/K$ be a connected and simply connected compact homogeneous space. Then $G/K$ is a generalized Wallach space if and only if one of the following types:

  1) $G/K$  is a direct product of three irreducible symmetric spaces of compact type ($[ijk]=0$ in this case).

  2) The group is simple and the pair $(\mathfrak{g},\mathfrak{k})$ is   one of the pairs in Table 1.

  3) $G=F\times F\times F\times F$ and $H=diag(F)\subset G$ for some connected simple connected compact simple Lie group $F$, with the following description on the Lie algebra level:
  $$
  (\mathfrak{g},\mathfrak{k})=(\mathfrak{f}\oplus \mathfrak{f}\oplus\mathfrak{f}\oplus \mathfrak{f}, diag(\mathfrak{f})=\{(X,X,X,X)\mid X \in f\},
  $$
  where $\mathfrak{f}$ is the Lie algebra of $F$, and (up to permutation) $\mathfrak{m}_1=\{(X,X,-X,-X)\mid X \in f\}$, $\mathfrak{m}_2=\{(X,-X,X,-X)\mid X \in f\}$,  $\mathfrak{m}_3=\{(X,-X,-X,X)\mid X \in f\}$.
  \end{theorem}

\medskip
\begin{center}
\begin{tabular}{|c|c|c|c|}
 \hline
         $\mathfrak{g}$  & $\mathfrak{k}$   & $\mathfrak{g}$  & $\mathfrak{k}$  \\
     \thickline
        $\mathfrak{so}(k+l+m)$ & $\mathfrak{so}(k)\oplus \mathfrak{so}(l)\oplus \mathfrak{so}(m)$ & $\mathfrak{e}_7$&$\mathfrak{so}(8)\oplus 3 \mathfrak{sp}(1)$   \\
    \thickline
       $\mathfrak{su}(k+l+m)$ & $\mathfrak{su}(k)\oplus su(l)\oplus \mathfrak{su}(m)$  &$\mathfrak{e}_7$ &$\mathfrak{su}(6)\oplus \mathfrak{sp}(1)\oplus \mathbb{R}$ \\
    \thickline
     $\mathfrak{sp}(k+l+m)$  &$\mathfrak{sp}(k)\oplus \mathfrak{sp}(l)\oplus \mathfrak{sp}(m)$  &$\mathfrak{e}_7$&$\mathfrak{so}(8)$ \\
   \thickline
        $\mathfrak{su}(2l), l\geq 2$  & $\mathfrak{u}(l)$ & $\mathfrak{e}_8$  &$\mathfrak{so}(12)\oplus 2 \mathfrak{sp}(1)$ \\
    \thickline
       $\mathfrak{so}(2l), l\geq 4$ &$\mathfrak{u}(l)\oplus \mathfrak{u}(l-1)$ &$\mathfrak{e}_8$ & $\mathfrak{so}(8)\oplus \mathfrak{so}(8)$ \\
    \thickline
       $\mathfrak{e}_6$& $\mathfrak{su}(4)\oplus 2 \mathfrak{sp}(1)\oplus \mathbb{R}$ &$\mathfrak{f}_4$&$\mathfrak{so}(5)\oplus 2 \mathfrak{sp}(1)$ \\
    \thickline
        $\mathfrak{e}_6$ & $\mathfrak{so}(8)\oplus \mathbb{R}^2$ &$\mathfrak{f}_4$ &$\mathfrak{so}(8)$ \\
     \thickline
        $\mathfrak{e}_6$&  $\mathfrak{sp}(3)\oplus \mathfrak{sp}(1)$&&\\
  \hline
 \end{tabular}
 \end{center}
 \smallskip
 \begin{center}
  Table 1. \ {\small The pairs $(\mathfrak{g}, \mathfrak{k})$ corresponding to generalized Wallach spaces $G/K$ with $G$ simple.}
\end{center}

\section{g.o.  generalized Wallach spaces  }

Let $(G/K, g)$ be a generalized Wallach space with reductive decomposition  $\mathfrak{g}=\mathfrak{k}\oplus\mathfrak{m}$, equipped  with a $G$-invariant metric corresponding to a scalar product of the form (\ref{7}).  Let $l=\mathrm{dim}\mathfrak{k}$ and $d_i=\mathrm{dim}(\mathfrak{m}_i)$  $(i =1,2,3)$, and let $\{e^{0}_{i}\}$ and $\{e^{j}_{s}\}$ be  orthogonal bases of  $\mathfrak{k}$ and $\mathfrak{m}_j$ respectively with respect to $B$, where
  $ 1\leq i \leq l$, $j=1,2,3$ and $1 \leq s \leq d_j$. For any $X \in \mathfrak{g}\setminus\{0\}$ we write
  \begin{equation}\label{11}
  X=\sum^{l}_{i=1}x_{i}e^{0}_{i}+\sum^{d_1}_{j=1}a_{j}e^{1}_{j}+\sum^{d_2}_{k=1}b_{k}e^{2}_{k}+\sum^{d_3}_{s=1}c_{s}e^{3}_{s},\quad x_i, a_j, b_k, c_s \in \mathbb{R}.
  \end{equation}

  Then by Lemma \ref{L1}  it follows that $X$  is a geodesic vector if and only if
  \begin{equation}\label{12}
  \langle [\sum^{l}_{i=1}x_{i}e^{0}_{i}+\sum^{d_1}_{j=1}a_{j}e^{1}_{j}+\sum^{d_2}_{k=1}b_{k}e^{2}_{k}+\sum^{d_3}_{s=1}c_{s}e^{3}_{s},Y]_{\mathfrak{m}},
  \sum^{d_1}_{j=1}a_{j}e^{1}_{j}+\sum^{d_2}_{k=1}b_{k}e^{2}_{k}+\sum^{d_3}_{s=1}a_{s}e^{3}_{s}\rangle=0,
  \end{equation}
  for all $Y\in \mathfrak{m}$.
  Hence we obtain  the following  system of $d_1+d_2+d_3$ equations
  \begin{equation}\label{13}
\begin{cases}
-\sum^{l}_{i=1}(\sum^{d_1}_{j=1}a_{j}B([e^{0}_{i},e^{1}_{1}]_{\mathfrak{m}},e^{1}_{j}))x_i\lambda_1=\sum^{d_2}_{k=1}\sum^{d_3}_{s=1}
b_kc_sB([e^{2}_{k},e^{1}_{1}]_{\mathfrak{m}},e^{3}_{s})(\lambda_3-\lambda_2)\\
\; \quad \quad \quad\quad\quad\quad\quad\quad \vdots \; \quad \quad \quad\quad\quad\quad\quad\quad \quad \quad \quad\quad\quad\quad\quad\quad \vdots\\
 -\sum^{l}_{i=1}(\sum^{d_1}_{j=1}a_{j}B([e^{0}_{i},e^{1}_{d_1}]_{\mathfrak{m}},e^{1}_{j}))x_i\lambda_1=\sum^{d_2}_{k=1}\sum^{d_3}_{s=1}
 b_kc_sB([e^{2}_{k},e^{1}_{d_1}]_{\mathfrak{m}},e^{3}_{s})(\lambda_3-\lambda_2)\\
 -\sum^{l}_{i=1}(\sum^{d_2}_{k=1}b_{k}B([e^{0}_{i},e^{2}_{1}]_{\mathfrak{m}},e^{2}_{k}))x_i\lambda_2=\sum^{d_1}_{j=1}\sum^{d_3}_{s=1}
 a_jc_sB([e^{1}_{j},e^{2}_{1}]_{\mathfrak{m}},e^{3}_{s})(\lambda_3-\lambda_1)\\
 \; \quad \quad \quad\quad\quad\quad\quad\quad \vdots \; \quad \quad \quad\quad\quad\quad\quad\quad \quad \quad \quad\quad\quad\quad\quad\quad \vdots\\
 -\sum^{l}_{i=1}(\sum^{d_2}_{k=1}b_{k}B([e^{0}_{i},e^{2}_{d_2}]_{\mathfrak{m}},e^{2}_{k}))x_i\lambda_2=\sum^{d_1}_{j=1}\sum^{d_3}_{s=1}
 a_jc_sB([e^{1}_{j},e^{2}_{d_2}]_{\mathfrak{m}},e^{3}_{s})(\lambda_3-\lambda_1)\\
  -\sum^{l}_{i=1}(\sum^{d_3}_{s=1}c_{s}B([e^{0}_{i},e^{3}_{1}]_{\mathfrak{m}},e^{3}_{s}))x_i\lambda_3=\sum^{d_1}_{j=1}\sum^{d_2}_{k=1}
 a_jb_kB([e^{1}_{j},e^{3}_{1}]_{\mathfrak{m}},e^{2}_{k})(\lambda_2-\lambda_1)\\
 \; \quad \quad \quad\quad\quad\quad\quad\quad \vdots \; \quad \quad \quad\quad\quad\quad\quad\quad \quad \quad \quad\quad\quad\quad\quad\quad \vdots\\
 -\sum^{l}_{i=1}B(\sum^{d_3}_{s=1}c_{s}B([e^{0}_{i},e^{3}_{d_3}]_{\mathfrak{m}},e^{3}_{s}))x_i\lambda_3=\sum^{d_1}_{j=1}\sum^{d_2}_{k=1}
 a_jb_kB([e^{1}_{j},e^{3}_{d_3}]_{\mathfrak{m}},e^{2}_{k})(\lambda_2-\lambda_1).\\
\end{cases}
\end{equation}

If we set
$$
A=-
\left (
\begin{array}{ccccccccccc}
\sum^{d_1}_{j=1}a_{j}B([e^{0}_{1},e^{1}
_{1}]_{\mathfrak{m}},e^{1}_{j})\lambda_1& \cdots & \sum^{d_1}_{j=1}a_{j}B([e^{0}_{l},e^{1}_{1}]_{\mathfrak{m}},e^{1}_{j})\lambda_1\\
\vdots & & \vdots\\
\sum^{d_1}_{j=1}a_{j}B([e^{0}_{1},e^{1}_{d_1}]_{\mathfrak{m}},e^{1}_{j})\lambda_1&\cdots &\sum^{d_1}_{j=1}a_{j}B([e^{0}_{l},e^{1}_{d_1}]_{\mathfrak{m}},e^{1}_{j})\lambda_1\\
\sum^{d_2}_{k=1}b_{k}B([e^{0}_{1},e^{2}_{1}]_{\mathfrak{m}},e^{2}_{k})\lambda_2 &\cdots &\sum^{d_2}_{k=1}b_{k}B([e^{0}_{l},e^{2}_{1}]_{\mathfrak{m}},e^{2}_{k})\lambda_2 \\
\vdots  & &\vdots\\
\sum^{d_2}_{k=1}b_{k}B([e^{0}_{1},e^{2}_{d_2}]_{\mathfrak{m}},e^{2}_{k})\lambda_2 & \cdots&\sum^{d_2}_{k=1}b_{k}B([e^{0}_{l},e^{2}_{d_2}]_{\mathfrak{m}},e^{2}_{k})\lambda_2\\
\sum^{d_3}_{s=1}c_{s}B([e^{0}_{1},e^{3}_{1}]_{\mathfrak{m}},e^{3}_{s})\lambda_3 &\cdots &\sum^{d_3}_{s=1}c_{s}B([e^{0}_{l},e^{3}_{1}]_{\mathfrak{m}},e^{3}_{s})\lambda_3\\
\vdots&& \vdots\\
\sum^{d_3}_{s=1}c_{s}B([e^{0}_{1},e^{3}_{d_3}]_{\mathfrak{m}},e^{3}_{s})\lambda_3& \cdots&\sum^{d_3}_{s=1}c_{s}B([e^{0}_{l},e^{3}_{d_3}]_{\mathfrak{m}},e^{3}_{s})\lambda_3\\
\end{array}
\right ),
$$
$$
\texttt{X}=
\left (
\begin{array}{ccccccccccc}
x_1&\\
 x_2&\\
\vdots\\
x_l&
\end{array}
\right ),
$$
and
$$
B=
\left (
\begin{array}{ccccccccccc}
\sum^{d_2}_{k=1}\sum^{d_3}_{s=1}b_kc_sB([e^{2}_{k},e^{1}_{1}]_{\mathfrak{m}},e^{3}_{s})(\lambda_3-\lambda_2)&\\
 \vdots\\
\sum^{d_2}_{k=1}\sum^{d_3}_{s=1}b_kc_sB([e^{2}_{k},e^{1}_{d_1}]_{\mathfrak{m}},e^{3}_{s})(\lambda_3-\lambda_2)&\\
\sum^{d_1}_{j=1}\sum^{d_3}_{s=1}a_jc_sB([e^{1}_{j},e^{2}_{1}]_{\mathfrak{m}},e^{3}_{s})(\lambda_3-\lambda_1)& \\
\vdots\\
\sum^{d_1}_{j=1}\sum^{d_3}_{s=1}a_jc_sB([e^{1}_{j},e^{2}_{d_2}]_{\mathfrak{m}},e^{3}_{s})(\lambda_3-\lambda_1)&\\
\sum^{d_1}_{j=1}\sum^{d_2}_{k=1}
 a_jb_kB([e^{1}_{j},e^{3}_{1}]_{\mathfrak{m}},e^{2}_{k})(\lambda_2-\lambda_1)&\\
 \vdots\\
\sum^{d_1}_{j=1}\sum^{d_2}_{k=1}
 a_jb_kB([e^{1}_{j},e^{3}_{d_3}]_{\mathfrak{m}},e^{2}_{k})(\lambda_2-\lambda_1)&\\
\end{array}
\right ),
$$
then System (\ref{13}) is equivalent to $A\texttt{X}=B$. Then we have the following:
\begin{proposition}\label{P2}
Let $(G/K, g)$ be a generalized  Wallach space, where the metric $g$ is determined by the scalar product (\ref{7}). Then for any  $a_j$ $(j=1,\dots,d_1)$, $b_k$ $(k=1,\dots,d_2)$ and $c_s $ $(s=1,\dots,d_3)$  not all equal to zero,  $(G/K, g)$ is  a g.o. space if and only if $\mathrm{rank}(A)=\mathrm{rank}(A,B)$.
\end{proposition}

 \begin{proof}
  By Propostion \ref{P1}  and Corollary \ref{C1} it follows  that $(G/K, g)$ is a g.o. space  if and only if for any $X\in \mathfrak{m}\setminus \{0\}$  there exists an  $a(X)$ such that the curve $\mathrm{exp}t(a(X)+X)\cdot o\;(t\in \mathbb{R})$ is a geodesic. This means that System (\ref{13}) of $d_1+d_2+d_3$ equations for the variables $x_i \;(i=1,\dots,l)$ has real solutions for any $a_j, b_k, c_s \;(j=1,\dots,d_1,k=1,\dots,d_2, s=1,\dots,d_3)$ not all equal to zero.  Then   system  $A\texttt{X}=B$ has real solutions if and only if $\mathrm{rank}(A)=\mathrm{rank}(A,B).$
\end{proof}

Next, we shall investigate which of the families of spaces listed in Theorem \ref{L1} are g.o. spaces.\\

\begin{theorem}\label{T3} Let $(G/K, g)$ be a generalized Wallach space as listed in Theorem \ref{L1}. Then

1)\ If $(G/K, g)$ is a space  of  type 1)  then this is a  g.o. space  for any $\mathrm{Ad}(K)$-invariant Riemannian metric.

2) \ If  $(G/K, g)$ is a space  of  type 2) or 3) then this  is a g.o. space if and only if $g$ is the standard metric.
\end{theorem}

\begin{proof} 1)  In this case it is $[ijk]=0$  for all $i, j, k$, therefore  $B([e^{i}_{\alpha}, e^{j}_{\beta}], e^{k}_{\gamma})=0$ for any $i,j,k$ and $\alpha, \beta, \gamma$. So for any $\Ad({K})$-invariant Riemannian metric $(\lambda_1, \lambda_2, \lambda_3)$ and  any non zero vector $X\in \mathfrak{m}\setminus \{0\}$, system (\ref{13}) for the variables $x_i \;(i=1,\dots,l)$ has  real solutions, hence $(G/K, g)$ is a g.o. space.

2) The only non zero $[ijk]$ is  $[123]$ (cf. \cite{{Nik1}}),  so $B([e^{1}_{\alpha},e^{2}_{\beta}],e^{3}_{\gamma}) \neq 0$ for some $\alpha, \beta, \gamma$. Let $(\lambda_1, \lambda_2,\lambda_3)$ be a g.o. metric,  that is  the generalized Wallach space is a g.o. space.
Then for any non zero vector $X \in \mathfrak{m}$,  there exists $a(X)$,  such that $X +a(X)$  is a geodesic vector.

 We assume  that  $X \in \mathfrak{m}_1 \oplus \mathfrak{m}_2$. Then   in  system (\ref{13}) all $c_1, c_2,\dots, c_{d_3}$ are  zero, hence it follows that

$$
A=-
\left (
\begin{array}{ccccccccccc}
\sum^{d_1}_{j=1}a_{j}B([e^{0}_{1},e^{1}_{1}]_{\mathfrak{m}},e^{1}_{j})\lambda_1& \cdots & \sum^{d_1}_{j=1}a_{j}B([e^{0}_{l},e^{1}_{1}]_{\mathfrak{m}},e^{1}_{j})\lambda_1\\
\vdots && \vdots\\
\sum^{d_1}_{j=1}a_{j}B([e^{0}_{1},e^{1}_{d_1}]_{\mathfrak{m}},e^{1}_{j})\lambda_1&\cdots &\sum^{d_1}_{j=1}a_{j}B([e^{0}_{l},e^{1}_{d_1}]_{\mathfrak{m}},e^{1}_{j})\lambda_1\\
\sum^{d_2}_{k=1}b_{k}B([e^{0}_{1},e^{2}_{1}]_{\mathfrak{m}},e^{2}_{k})\lambda_2 &\cdots &\sum^{d_2}_{k=1}b_{k}B([e^{0}_{l},e^{2}_{1}]_{\mathfrak{m}},e^{2}_{k})\lambda_2 \\
\vdots  & &\vdots\\
\sum^{d_2}_{k=1}b_{k}B([e^{0}_{1},e^{2}_{d_2}]_{\mathfrak{m}},e^{2}_{k})\lambda_2 & \cdots&\sum^{d_2}_{k=1}b_{k}B([e^{0}_{l},e^{2}_{d_2}]_{\mathfrak{m}},e^{2}_{k})\lambda_2\\
0 &\cdots &0\\
\vdots&& \vdots\\
0& \cdots&0\\
\end{array}
\right ),
$$
\\
$$
B=
\left (
\begin{array}{ccccccccccc}
0&\\
 \vdots\\
0&\\
0\\
\vdots\\
0&\\
\sum^{d_1}_{j=1}\sum^{d_2}_{k=1}
 a_jb_kB([e^{1}_{j},e^{3}_{1}]_{\mathfrak{m}},e^{2}_{k})(\lambda_2-\lambda_1)&\\
 \vdots\\
\sum^{d_1}_{j=1}\sum^{d_2}_{k=1}
 a_jb_kB([e^{1}_{j},e^{3}_{d_3}]_{\mathfrak{m}},e^{2}_{k})(\lambda_2-\lambda_1)&\\
\end{array}
\right ).
$$\\
 Then it is $\mathrm{rank}(A)=\mathrm{rank}(A,B)$   for any $a_j, b_k\;(j=1,\dots,d_1,k=1,\dots,d_2)$  not all equal to zero    if and only if $\lambda_1 =\lambda_2$.

 We now assume that  $X \in \mathfrak{m}_1 \oplus \mathfrak{m}_3$.  Then  in system (\ref{13}) all $b_1, b_2,\dots, b_{d_3}$ are zero, hence we obtain that

 $$
A=-
\left (
\begin{array}{ccccccccccc}
\sum^{d_1}_{j=1}a_{j}B([e^{0}_{1},e^{1}_{1}]_{\mathfrak{m}},e^{1}_{j})\lambda_1& \cdots & \sum^{d_1}_{j=1}a_{j}B([e^{0}_{l},e^{1}_{1}]_{\mathfrak{m}},e^{1}_{j})\lambda_1\\
\vdots && \vdots\\
\sum^{d_1}_{j=1}a_{j}B([e^{0}_{1},e^{1}_{d_1}]_{\mathfrak{m}},e^{1}_{j})\lambda_1&\cdots &\sum^{d_1}_{j=1}a_{j}B([e^{0}_{l},e^{1}_{d_1}]_{\mathfrak{m}},e^{1}_{j})\lambda_1\\
0 &\cdots &0 \\
\vdots  & &\vdots\\
0& \cdots&0\\
\sum^{d_3}_{s=1}c_{s}B([e^{0}_{1},e^{3}_{1}]_{\mathfrak{m}},e^{3}_{s})\lambda_3 &\cdots &\sum^{d_3}_{s=1}c_{s}B([e^{0}_{l},e^{3}_{1}]_{\mathfrak{m}},e^{3}_{s})\lambda_3\\
\vdots&& \vdots\\
\sum^{d_3}_{s=1}c_{s}B([e^{0}_{1},e^{3}_{d_3}]_{\mathfrak{m}},e^{3}_{s})\lambda_3& \cdots&\sum^{d_3}_{s=1}c_{s}B([e^{0}_{l},e^{3}_{d_3}]_{\mathfrak{m}},e^{3}_{s})\lambda_3\\
\end{array}
\right ),
$$
\\
$$
B=
\left (
\begin{array}{ccccccccccc}
0&\\
 \vdots\\
0&\\
\sum^{d_1}_{j=1}\sum^{d_3}_{s=1}a_jc_sB([e^{1}_{j},e^{2}_{1}]_{\mathfrak{m}},e^{3}_{s})(\lambda_3-\lambda_1)& \\
\vdots\\
\sum^{d_1}_{j=1}\sum^{d_3}_{s=1}a_jc_sB([e^{1}_{j},e^{2}_{d_2}]_{\mathfrak{m}},e^{3}_{s})(\lambda_3-\lambda_1)&\\
0&\\
 \vdots\\
0&\\
\end{array}
\right ).
$$

Therefore,  $\mathrm{rank}(A)=\mathrm{rank}(A,B)$   for any $a_j, c_s \; (j=1,\dots,d_1,s=1,\dots,d_3, )$  not all equal to zero,  if and only if $\lambda_1=\lambda_3$.
We conclude that the generalized Wallach spaces of types 2), 3) in Theorem \ref{L1} are g.o. spaces, if and only  if  the  metric  $(\lambda_1, \lambda_2, \lambda_3)$  satisfies  $\lambda_1=\lambda_2=\lambda_3$.
\end{proof}

\section{Homogeneous geodesics in  generalized Wallach spaces}
Let $(G/K, g)$ be a generalized Wallach space with $B$-orthogonal decomposition  $\mathfrak{g}=\mathfrak{k}\oplus\mathfrak{m}$,  where $ \mathfrak{m}=\mathfrak{m}_1\oplus\mathfrak{m}_2\oplus\mathfrak{m}_3$ and $B=-$Killing form on $\mathfrak{g}$.
Let $\{e^{0}_{i}\}$ and $\{e^{j}_{s}\}$ be the orthogonal bases of $\mathfrak{k}$ and $\mathfrak{m}_j$ respectively with respect to $B$
  ($1\leq i \leq l$, $j=1,2,3$, $1 \leq s \leq d_j$).  For any $X \in \mathfrak{g}\setminus\{0\}$ we write
  $$
  X=\sum^{l}_{i=1}x_{i}e^{0}_{i}+\sum^{d_1}_{j=1}a_{j}e^{1}_{j}+\sum^{d_2}_{k=1}
  b_{k}e^{2}_{k}+\sum^{d_3}_{s=1}c_{s}e^{3}_{s},\qquad
 x_i, a_j, b_k, c_s \in \mathbb{R}.
  $$

In order to  find all  homogeneous geodesics in $G/K$, it  suffices  to find all the real solutions of  system (\ref{13}), of $d_1+d_2+d_3$ quadratic equations for the variables $x_i, a_j, b_k, c_s$, which are  not all equal to zero.

 By Theorem \ref{T3} we only consider homogeneous geodesics in generalized Wallach spaces types 2) and 3) given in Theorem \ref{T1} for the metric $(\lambda_{1}, \lambda_{2}, \lambda_{3})$, where  at least two of  $\lambda_{1}, \lambda_{2}, \lambda_{3}$ are  different.  Then the geodesic vectors correspond to those  solutions of  system (\ref{13}) for the variables $x_i, a_j, b_k, c_s$, which are  not all equal to zero. However, for many generalized Wallach spaces it is difficult to find all the real solutions of system (\ref{13}).

\medskip
Next, we will give two examples of generalized Wallach spaces and give all the homogeneous geodesics for any given metric.
\begin{example}\label{E1}
  We consider the generalized Wallach space $SU(2)/\{e\}$.
  \end {example}

Let $\{\sqrt{-1}h_{\alpha}, \frac{A_{\alpha}}{\sqrt{2}}, \frac{B_{\alpha}}{\sqrt{2}} \}$ be an orthogonal basis of $\mathfrak{su}(2)$ with respect to $B$,  where $\alpha$ denotes a simple root of the Lie algebra  $\mathfrak{su}(2)$, and $A_{\alpha}=E_{\alpha}-E_{-\alpha}, \ B_{\alpha}=\sqrt{-1}(E_{\alpha}+E_{-\alpha})$.  Here $\{ E_{\alpha}  \}$ denotes the Weyl basis of  $\mathfrak{su}(2)$.
We set  $X_{\alpha}=\frac{A_{\alpha}}{\sqrt{2}}, Y_{\alpha}=\frac{B_{\alpha}}{\sqrt{2}}$. Then we have that
$$
 \quad [\sqrt{-1}h_{\alpha}, X_{\alpha}]=\alpha(h_{\alpha})Y_{\alpha}=Y_{\alpha},
$$
$$
\quad \quad \quad [\sqrt{-1}h_{\alpha}, Y_{\alpha}]=-\alpha(h_{\alpha})X_{\alpha}=-X_{\alpha},
$$
$$
[X_{\alpha}, Y_{\alpha}]=\sqrt{-1}h_{\alpha}.
$$

The Lie algebra  $\mathfrak{su}(2)$ has an orthogonal decomposition  $\mathfrak{su}(2)=\mathfrak{m}_1\oplus\mathfrak{m}_2\oplus\mathfrak{m}_3$, where $\mathfrak{m}_1, \mathfrak{m}_2, \mathfrak{m}_3$ are spanned by $\sqrt{-1}h_{\alpha},  X_{\alpha} $ and $Y_{\alpha}$ respectively.  For any $X \in \mathfrak{su}(2) $ we write
$$X=a\sqrt{-1}h_{\alpha}+bX_{\alpha}+cY_{\alpha}, \qquad a, b, c \in \mathbb{R}.
$$

 We will  find all  geodesic vectors of   $SU(2)/\{e\}$ for a given metric $(\lambda _1, \lambda _2, \lambda _3)$. System (\ref{13}) for $SU(2)/\{e\}$ is

\begin{equation}\label{18}
\begin{cases}
bc(\lambda_3-\lambda_2)=0\\
ac(\lambda_3-\lambda_1)=0\\
ab(\lambda_2-\lambda_1)=0.\\
\end{cases}
\end{equation}

\noindent
 {\bf Case 1.} $\lambda_1=\lambda_2\neq \lambda_3$.
The solutions of the system (\ref{18}) are  $c=0$ or $c\neq 0, a=b=0$, so  the geodesic vectors are $X=a\sqrt{-1}h_{\alpha}+bX_{\alpha}$ and $X=cY_{\alpha}$.

\noindent
 {\bf Case 2.}  $\lambda_1=\lambda_3\neq \lambda_2$.
 The solutions of the system (\ref{18}) are  $b=0$ or $b\neq 0, a=c=0$,  so  the geodesic vectors are $X=a\sqrt{-1}h_{\alpha}+cY_{\alpha}$  and $X=bX_{\alpha}$.

\noindent
 {\bf Case 3.}  $\lambda_2=\lambda_3\neq \lambda_1$.
 The solutions of the system (\ref{18}) are  $a=0$ or $a\neq 0, b=c=0$,  so  the geodesic vectors  are $X=bX_{\alpha}+cY_{\alpha}$ and  $X=a\sqrt{-1}h_{\alpha}$.

\noindent
 {\bf Case 4.}  $\lambda_1, \lambda_2, \lambda_3$ are distinct.
Then system (\ref{18}) reduces to $ab=ac=bc=0$, whose solutions are $a=b=0$ or $a=c=0$ or $b=c=0$. In this case any vector $X \in \mathfrak{m}_i\backslash\{0\}$ $(i=1, 2, 3)$ is a geodesic vector.

\smallskip
Therefore we obtain the following theorem, which   recovers a result on R.A. Marinosci \cite[p. 266]{Mar}

\begin{proposition}\label{T4} For the generalized Wallach space ${SU}(2)/\{e\}$ the only geodesic vectors for a given metric $(\lambda_1, \lambda_2, \lambda_3)$ are the following:

\noindent
1)\ If $\lambda_i=\lambda_j \neq \lambda_k \;(i, j, k \in \{1, 2, 3\})$, then any vector $X\in \mathfrak{m}_k \backslash\{0\}$ or $X\in (\mathfrak{m}_i\oplus \mathfrak{m}_j)\backslash \{0\}$.

\noindent
2)\ If $\lambda_1, \lambda_2, \lambda_3$ are distinct, then any vector $X\in \mathfrak{m}_1\cup\mathfrak{m}_2\cup\mathfrak{m}_3$.
\end{proposition}

\smallskip
\begin{example}\label{E2}
 We consider the generalized Wallach space $SO(n)/SO(n-2)$, ($n\ge 4$).
 \end{example}

This is the Stiefel manifold  of orthonormal $2$-frames in $\mathbb{R}^n$.
Let $\mathfrak{so}(n)$ and $\mathfrak{so}(n-2)$ be the Lie algebras of $SO(n)$ and $SO(n-2)$ respectively. Let $E_{ab}$ denote the $n\times n$ matrix with
 1 in the (ab)-entry and 0 elsewhere.  If $e_{ab}=E_{ab}-E_{ba}$, then the set $\mathcal{B}=\{e_{ij}=E_{ij}-E_{ji}: 1\leq i < j\leq n \}$ is a  $B$-orthogonal basis of $\mathfrak{so}(n)$.
 The multiplication table of the elements in $\mathcal{B}$ is given as follows:

 \begin{lemma} If all four indices are distinct, then the Lie brackets in $\mathcal{B}$ are zero.  Otherwise, it is
 $[e_{ij}, e_{jk}]=e_{ik}$, where $i, j, k$ are distinct.
 \end{lemma}

 Let  $\mathfrak{so}(n)=\mathfrak{k}\oplus \mathfrak{m}_1\oplus\mathfrak{m}_2\oplus\mathfrak{m}_3$ be an orthogonal decomposition of $\mathfrak{so}(n)$ with respect to
 $B$, where $\mathfrak{k}={\mathrm{span}}_{\mathbb{R}}\{e_{ij}: 3\leq i < j\leq n\}$,  $\mathfrak{m}_1=\mathrm{span}_{\mathbb{R}}\{e_{12}\}$, $\mathfrak{m}_2=\mathrm{span}_{\mathbb{R}}\{e_{1j}:3 \leq j\leq n\}$, and  $\mathfrak{m}_3=\mathrm{span}_{\mathbb{R}}\{e_{2j}:3 \leq j\leq n\}$.  For any $X \in \mathfrak{so}(n) $ we write
 $$
 X=\sum_{3\leq i < j\leq n}a_{ij}e_{ij}+a_{12}e_{12}+\sum_{3 \leq j\leq n}a_{1j}e_{1j}+
 \sum_{3 \leq j\leq n}a_{2j}e_{2j}, \qquad a_{ij}\in \mathbb{R}.
 $$
Then system (\ref{13}) for   $SO(n)/SO(n-2)$ takes the form

 \begin{equation}\label{20}
\begin{cases}
(a_{13}a_{23}+a_{14}a_{24}+\cdots+a_{1n}a_{2n})(\lambda_3-\lambda_2)=0\\
-(a_{14}a_{34}+a_{15}a_{35}+\cdots+a_{1, n}a_{3n})\lambda_2=a_{12}a_{23}(\lambda_3-\lambda_1)\\
 (a_{13}a_{34}-a_{15}a_{45}-a_{16}a_{46}-\cdots-a_{1n}a_{4n})\lambda_2=a_{12}a_{24}(\lambda_3-\lambda_1)\\
 (a_{13}a_{35}+a_{14}a_{45}-a_{16}a_{56}-a_{17}a_{57}-\cdots-a_{1n}a_{5n})\lambda_2=a_{12}a_{25}(\lambda_3-\lambda_1)\\
 \quad\quad\quad\quad\quad\quad\quad\quad\quad\quad\quad\quad\quad\quad\quad\vdots\\
(a_{13}a_{3, n-1}+\cdots+ a_{1, n-2}a_{n-2, n-1}-a_{1n}a_{n-1, n})\lambda_2=a_{12}a_{2, n-1}(\lambda_3-\lambda_1)\\
(a_{13}a_{3n}+a_{14}a_{4,n}+\cdots+ a_{1, n-1}a_{n-1, n})\lambda_2=a_{12}a_{2n}(\lambda_3-\lambda_1)\\
-(a_{24}a_{34}+a_{25}a_{35}+\cdots+a_{2n}a_{3n})\lambda_3=a_{12}a_{13}(\lambda_1-\lambda_2)\\
 (a_{23}a_{34}-a_{25}a_{45}-a_{26}a_{46}-\cdots-a_{2n}a_{4n})\lambda_3=a_{12}a_{14}(\lambda_1-\lambda_2)\\
 (a_{23}a_{35}+a_{24}a_{45}-a_{26}a_{56}-a_{27}a_{57}-\cdots-a_{2n}a_{5n})\lambda_3=a_{12}a_{15}(\lambda_1-\lambda_2)\\
 \quad\quad\quad\quad\quad\quad\quad\quad\quad\quad\quad\quad\quad\quad\quad\vdots\\
(a_{23}a_{3, n-1}+\cdots+ a_{2, n-2}a_{n-2, n-1}-a_{2n}a_{n-1, n})\lambda_3=a_{12}a_{1, n-1}(\lambda_1-\lambda_2)\\
(a_{23}a_{3n}+a_{24}a_{4n}+\cdots+ a_{2, n-1}a_{n-1, n})\lambda_3=a_{12}a_{1, n}(\lambda_1-\lambda_2).\\
\end{cases}
\end{equation}

As the above system is difficult to handle, we restrict to the  Stiefel manifold
  $SO(4)/SO(2)$ and look for geodesics
  $X=a_{34}e_{34}+a_{12}e_{12}+a_{13}e_{13}+a_{14}e_{14}+a_{23}e_{23}+a_{24}e_{24}.$
  Then the above system simplifies to

 \begin{equation}\label{20}
\begin{cases}
(a_{13}a_{23}+a_{14}a_{24})(\lambda_3-\lambda_2)=0\\
a_{34}a_{14}\lambda_2=-a_{12}a_{23}(\lambda_3-\lambda_1)\\
 a_{34}a_{13}\lambda_2=a_{12}a_{24}(\lambda_3-\lambda_1)\\
a_{34}a_{24}\lambda_3=a_{13}a_{12}(\lambda_2-\lambda_1)\\
a_{34}a_{23}\lambda_3=-a_{14}a_{12}(\lambda_2-\lambda_1).\\
\end{cases}
\end{equation}

 In order  to find all  geodesic vectors in the generalized Wallach space  $SO(4)/SO(2)$ for a given metric, we should find all the non zero real solutions of  the system (\ref{20}).

\smallskip
\noindent
{\bf Case 1.} $\lambda_1=\lambda_2\neq \lambda_3$.
Then system  (\ref{20}) reduces to

\begin{equation}\label{21}
\begin{cases}
a_{13}a_{23}+a_{14}a_{24}=0\\
a_{34}a_{14}\lambda_1=-a_{12}a_{23}(\lambda_3-\lambda_1)\\
 a_{34}a_{13}\lambda_1=a_{12}a_{24}(\lambda_3-\lambda_1).\\
\end{cases}
\end{equation}

If $a_{34}=0$  and  $a_{12}=0$ then the geodesic vectors  are
$X=a_{13}e_{13}+a_{14}e_{14}+a_{23}e_{23}+a_{24}e_{24}$  with $a_{13}a_{23}+a_{14}a_{24}=0$.

 If $a_{34}=0$  and $a_{12}\neq 0$ we get $a_{23}=a_{24}=0$, so the   geodesic vectors  are $X=a_{12}e_{12}+a_{13}e_{13}+a_{14}e_{14}$.

If $a_{34}\neq 0$ and $a_{12}=0$  we get $a_{13}=a_{14}=0$, so the  geodesic vectors  are $X=a_{34}e_{34}+a_{23}e_{23}+a_{24}e_{24}$.

If $a_{34}\neq 0$, $a_{12}\neq0$ and  $a_{23}+a_{24}=0$  we have $a_{13}=a_{14}$,  geodesic vectors  are  $X=a_{34}e_{34}+a_{12}e_{12}+a_{13}e_{13}+a_{13}e_{14}+a_{23}e_{23}-a_{23}e_{24}$.

If $a_{34}\neq 0$, $a_{12}\neq0$ and $a_{23}+a_{24}\neq 0$ we have  $a_{13}\neq a_{14}$.
If  $a_{13}= -a_{14}$ we have  $a_{23}=a_{24}$,  so  geodesic vectors  are $X=a_{34}e_{34}+a_{12}e_{12}+a_{13}e_{13}-a_{13}e_{14}+a_{23}e_{23}+a_{23}e_{24}$.  If  $a_{13}\neq -a_{14}$  we have  $a_{23}\neq a_{24}$,  so geodesic vectors  are $X=a_{34}e_{34}+a_{12}e_{12}+a_{13}e_{13}+a_{14}e_{14}+a_{23}e_{23}+a_{24}e_{24}$.

\smallskip
\noindent
{\bf Case 2.} $\lambda_1=\lambda_3\neq \lambda_2$.
Then system  (\ref{20})  reduces to

\begin{equation}\label{23}
\begin{cases}
a_{13}a_{23}+a_{14}a_{24}=0\\
a_{34}a_{24}\lambda_1=a_{13}a_{12}(\lambda_2-\lambda_1)\\
 a_{34}a_{23}\lambda_1=-a_{12}a_{14}(\lambda_2-\lambda_1).\\
\end{cases}
\end{equation}

 If $a_{34}=0$ and $a_{12}=0$,   geodesic vectors  are  $X=a_{13}e_{13}+a_{14}e_{14}+a_{23}e_{23}+a_{24}e_{24}$  with $a_{13}a_{23}+a_{14}a_{24}=0$.

  If $a_{34}=0$  and  $a_{12}\neq0$  we get $a_{13}=a_{14}=0$,  geodesic vectors  are $X=a_{12}e_{12}+a_{23}e_{23}+a_{24}e_{24}$.

 If $a_{34}\neq0$  and  $a_{12}=0$ we get $a_{23}=a_{24}=0$,  geodesic vectors  are $X=a_{34}e_{34}+a_{13}e_{13}+a_{14}e_{14}$.

If $a_{34}\neq 0$, $a_{12}\neq0$ and  $a_{23}+a_{24}=0$  we have $a_{13}=a_{14}$,  geodesic vectors  are  $X=a_{34}e_{34}+a_{12}e_{12}+a_{13}e_{13}+a_{13}e_{14}+a_{23}e_{23}-a_{23}e_{24}$.

If $a_{34}\neq 0$, $a_{12}\neq0$  and $a_{23}+a_{24}\neq 0$ we have  $a_{13}\neq a_{14}$.
If  $a_{13}= -a_{14}$ we have  $a_{23}=a_{24}$,  so geodesic vectors  are $X=a_{34}e_{34}+a_{12}e_{12}+a_{13}e_{13}-a_{13}e_{14}+a_{23}e_{23}+a_{23}e_{24}$.
If  $a_{13}\neq -a_{14}$  we have  $a_{23}\neq a_{24}$, so geodesic vectors  are $X=a_{34}e_{34}+a_{12}e_{12}+a_{13}e_{13}+a_{14}e_{14}+a_{23}e_{23}+a_{24}e_{24}$.

\smallskip
\noindent
{\bf Case 3.} $\lambda_3=\lambda_2\neq \lambda_1$.
Then system  (\ref{13}) reduces to

\begin{equation}\label{26}
\begin{cases}
a_{34}a_{14}\lambda_2=-a_{12}a_{23}(\lambda_2-\lambda_1)\\
 a_{34}a_{13}\lambda_2=a_{12}a_{24}(\lambda_2-\lambda_1)\\
a_{34}a_{24}\lambda_2=a_{13}a_{12}(\lambda_2-\lambda_1)\\
a_{34}a_{23}\lambda_2=-a_{14}a_{12}(\lambda_2-\lambda_1).\\
\end{cases}
\end{equation}

If $a_{34}=0$ and $a_{12}=0$,   geodesic vectors  are $X=a_{13}e_{13}+a_{14}e_{14}+a_{23}e_{23}+a_{24}e_{24}$.

If $a_{34}=0$ and $a_{12}\neq 0$ we have $a_{13}=a_{14}=a_{23}=a_{24}=0$,  geodesic vectors  are $X=a_{12}e_{12}$.

If $a_{34}\neq 0$ and  $a_{12}=0$ we have $a_{13}=a_{14}=a_{23}=a_{24}=0$,  geodesic vectors  are $X=a_{34}e_{34}$.

If $a_{34}\neq 0$, $a_{12}\neq0$ and  $a_{23}+a_{24}=0$  we have $a_{13}=a_{14}$,  geodesic vectors  are  $X=a_{34}e_{34}+a_{12}e_{12}+a_{13}e_{13}+a_{13}e_{14}+a_{23}e_{23}-a_{23}e_{24}$.

If $a_{34}\neq 0$ , $a_{12}\neq0$ and $a_{23}+a_{24}\neq 0$ we have  $a_{13}\neq a_{14}$.
If  $a_{13}= -a_{14}$ we have  $a_{23}=a_{24}$,  so geodesic vectors  are $X=a_{34}e_{34}+a_{12}e_{12}+a_{13}e_{13}-a_{13}e_{14}+a_{23}e_{23}+a_{23}e_{24}$.
If  $a_{13}\neq -a_{14}$  we have  $a_{23}\neq a_{24}$, so geodesic vectors  are $X=a_{34}e_{34}+a_{12}e_{12}+a_{13}e_{13}+a_{14}e_{14}+a_{23}e_{23}+a_{24}e_{24}$.

\smallskip
\noindent
{\bf Case 4.} $\lambda_1, \lambda_2, \lambda_3$ are all different.

If $a_{34}=0$ and $a_{12}=0$,  geodesic vectors  are $X=a_{13}e_{13}+a_{14}e_{14}+a_{23}e_{23}+a_{24}e_{24}$ with $a_{13}a_{23}+a_{14}a_{24}=0$.

If $a_{34}=0$ and  $a_{12}\neq0$ we have $a_{13}= a_{14}=a_{23}=a_{24}=0$,     geodesic vectors  are  $X=a_{12}e_{12}$.

If $a_{34}\neq0$ and $a_{12}=0$ we have  $a_{13}=a_{14}=a_{23}=a_{24}=0$,   geodesic vectors  are $X=a_{34}e_{34}$.

If $a_{34}\neq 0$, $a_{12}\neq0$ and  $a_{23}+a_{24}=0$  we have $a_{13}=a_{14}$,  geodesic vectors  are  $X=a_{34}e_{34}+a_{12}e_{12}+a_{13}e_{13}+a_{13}e_{14}+a_{23}e_{23}-a_{23}e_{24}$.

If $a_{34}\neq 0$ , $a_{12}\neq0$ and $a_{23}+a_{24}\neq 0$ we have  $a_{13}\neq a_{14}$.
If  $a_{13}= -a_{14}$ we have  $a_{23}=a_{24}$,  so geodesic vectors  are $X=a_{34}e_{34}+a_{12}e_{12}+a_{13}e_{13}-a_{13}e_{14}+a_{23}e_{23}+a_{23}e_{24}$.
If  $a_{13}\neq -a_{14}$  we have  $a_{23}\neq a_{24}$, so geodesic vectors  are $X=a_{34}e_{34}+a_{12}e_{12}+a_{13}e_{13}+a_{14}e_{14}+a_{23}e_{23}+a_{24}e_{24}$.

\end{document}